\documentclass[a4paper, 12pt]{article}

\usepackage[cp1251]{inputenc}
\usepackage[T2A]{fontenc}

\usepackage[russian,english]{babel}
\usepackage{amsmath, amsfonts, amssymb, amsthm}

\usepackage{geometry}

\usepackage[pdftex]{graphicx}
\graphicspath{{pictures/}}

\usepackage{pgf,tikz}
\usetikzlibrary{arrows}

\theoremstyle{theorem}

\newtheorem*{theoremA}{Theorem A}
\newtheorem*{theoremB}{Theorem B}
\newtheorem*{theoremC}{Theorem C}
\newtheorem*{theoremD}{Theorem D}
\newtheorem*{theoremE}{Theorem E}
\newtheorem*{theoremF}{Theorem F}
\newtheorem*{theoremMain}{Theorem}

\theoremstyle{definition}

\numberwithin{equation}{section}

\usepackage{caption}
\DeclareCaptionLabelSeparator{dot}{. }
\title{Reverse isoperimetric inequality in two-dimensional Alexandrov spaces}

\author{Alexander A. Borisenko\footnote{ National Academy of Sciences of Ukraine;  E-mail: \texttt{aborisenk@gmail.com}}}

\date{}

\begin{document}

\maketitle

\begin{abstract}
We prove a reverse isoperimetric inequality for domains homeomorphic to a disc with the boundary of curvature bounded below lying in two-dimensional Alexandrov spaces of curvature $\geqslant c$. We also study the equality case.

\textbf{Keywords:} Alexandrov metric spaces;  isoperimetric inequality; $\lambda$-convex curve
\end{abstract}

Well-known isoperimetric inequality for the Euclidean plane states that the area $F$ and the length $L$ of the boundary of any plane domain with a rectifiable boundary satisfy the inequality
$$
L^2 - 4 \pi F \geqslant 0,
$$
and equality is attained only for a circle~\cite{Bla56}.

On the planes of constant curvature there is a similar theorem, and the following inequality holds~\cite{Ber05}:
$$
L^2 - 4\pi F + c F^2 \geqslant 0,
$$ 
where $c$ is the curvature of the plane.

In two-dimensional manifolds of bounded curvature for domains homeomorphic to a disc A.\,D.\,Alexandrov proved~\cite{Al45} the inequality
$$
F \leqslant \frac{L^2}{2 (2\pi - \omega^+)},
$$
where $\omega^+$ is a positive curvature of a domain. In the inequality above equality holds only if the domain is isometric to a lateral surface of a right circular cone with curvature $\omega^+ < 2\pi$ at the vertex.

The isoperimetric inequality for domains with a compact closure and bounded by a finite number of rectifiable curves in two-dimensional manifolds of bounded curvature was proved in~\cite{Ion69}.  

But if we don't put any conditions on the boundary curves, then the areas of the enclosed domains can be arbitrary close to zero, and the perimeters of these curves can be arbitrary large. If we assume that the boundary curve is $\lambda$-convex, then, given the perimeter of it, the area of the domain is bounded from below by some constant. We prove the following main theorem.

 \begin{theoremMain}
 \label{CBBthm}
Let $G$ be a domain homeomorphic to a disc and lying in a two-dimensional Alexandrov space of curvature (in the sense of Alexandrov) $\geqslant c$. If the boundary curve $\gamma$ of $G$ is $\lambda$-convex, and the perimeter of $\gamma$ is equal to $L$, then the area $F$ of the domain $G$ satisfies
\begin{enumerate}
\item
for $c = 0$,
\begin{equation}
\label{cbbeq1}
F \geqslant \frac{L}{2\lambda} - \frac{1}{\lambda^2} \sin \frac{L\lambda}{2};
\end{equation}

\item
for $c = k^2$,
\begin{equation}
\label{cbbeq2}
F \geqslant \frac{4}{k^2} \arctan\left(\frac{\lambda}{\sqrt{\lambda^2+k^2}}\tan\left(\frac{\sqrt{\lambda^2+k^2}}{4}L\right)\right)-\frac{L\lambda}{k^2};
\end{equation}

\item
\begin{enumerate}
\item
for $c = -k^2$, and $\lambda > k$,
\begin{equation}
\label{cbbeq3}
F \geqslant \frac{L \lambda }{k^2} - \frac{4}{k^2} \arctan\left(\frac{\lambda}{\sqrt{\lambda^2 - k^2}}\tan\left(\frac{\sqrt{\lambda^2 - k^2}}{4}L\right)\right);
\end{equation}

\item
for $c = -k^2$, and $\lambda = k$,
\begin{equation}
\label{cbbeq4}
A \geqslant   \frac{L}{k} - \frac{4}{k^2} \arctan\frac{k L}{4};
\end{equation}

\item
for $c = -k^2$, and $\lambda < k$,
\begin{equation}
\label{cbbeq5}
F \geqslant \frac{L \lambda}{k^2} - \frac{4}{k^2} \arctan\left(\frac{\lambda}{\sqrt{k^2 - \lambda^2}}\tanh \left(\frac{\sqrt{k^2 - \lambda^2}}{4}L\right)\right).
\end{equation}
\end{enumerate}
\end{enumerate}
In all inequalities~(\ref{cbbeq1})~-- (\ref{cbbeq5}) equality is attained if and only if the domain $G$ is a $\lambda$-convex lune of length $L$ lying on the plane of constant curvature equal to $c$.
\end{theoremMain} 

By $\lambda$-convex lune we understand a convex domain bounded by two arcs of constant geodesic curvature and of length $L/2$. In particular:
\begin{enumerate}
\item[1)]
For $c = 0$ the domain is bounded by two circular arcs of radius $1/\lambda$. In this case the perimeter $L \leqslant 2\pi/\lambda$.

\item[2)]
For $c = k^2 > 0$ the domain is also bounded by two circular arcs of curvature equal to $\lambda$, and the perimeter satisfies $L \leqslant 2\pi/\sqrt{\lambda^2 + k^2}$.

\item[3)]
\begin{enumerate}
\item[(a)]
For $c = -k^2 < 0$ and $\lambda > k$ the lune is also bounded by two arcs of a circle of curvature equal to $\lambda$, and the perimeter satisfies $L \leqslant 2\pi / \sqrt{\lambda^2 - k^2}$;

\item[(b)]
for $c = -k^2 < 0$, $\lambda = k$ the domain is bounded by arcs of horocycles, and the perimeter of the domain can be arbitrary;

\item[(c)]
for $c = -k^2 < 0$, $0 < \lambda < k$ the domain is bounded by two arcs of equidistants, the perimeter can also be arbitrary.
\end{enumerate}
\end{enumerate}

For domains in two-dimensional simply-connected spaces of constant curvature equal to $c$ the main theorem was proved when $c = 0$ in~\cite{BorDr14}, when $c = k^2$ in~\cite{BorDr15_1}, and when $c = -k^2$ in~\cite{Dr14}.

For $J$-holomorphic curves some variant of a reverse isoperimetric inequality was proved in~\cite{GrS14}. In particular, it's been shown that the length of the boundary of a $J$-holomorphic curve with Lagrangian boundary conditions is dominated be a constant times its area.


\section{Alexandrov spaces}

Let $R$ be a metric space. Recall that a \textit{curve} in the metric space $R$ is a continuous image of a segment. If a curve $\gamma$ is given by a map $x = x(t)$, $a \leqslant t \leqslant b$, then the \textit{length} of $\gamma$ is defined as
$$
L_\gamma = \sup \sum \rho(x(t_{k-1}), x(t_k)), \quad a \leqslant t_1 \leqslant t_2 \leqslant \ldots \leqslant b,
$$ 
where $\rho(x(t_{k-1}), x(t_k))$ is the distance between the points $x(t_{k-1})$ and $x(t_k)$ in the space $R$, and the supremum is taken over all finite sub-divisions by the points $t_k$ of the segment $[a,b]$. 

Suppose that any two points $P$ and $Q$ of the space $R$ can be joined by a rectifiable curve. Then one can define a distance between $P$ and $Q$ in $R$ by taken the infimum over all lengths of curves joining them. A metric space in which the distance is given is such a way is called a \textit{space with inner metric}.

A curve $\gamma$ in the space with inner metric is called a \textit{shortest path} of a \textit{segment}, if its length is equal to the distance between its end points.

From now on we consider a two-dimensional space with inner metric.

We will call a \textit{triangle} a closed domain homeomorphic to a disc and whose boundary is composed of three shortest paths (sides of the triangle); between the sides one can define a lower angle.

We will say that a manifold posses a metric of curvature $\geqslant c$, if this metric is inner and for any sufficiently small triangle the sum of lower angles between its sides is not less than the sum of angles in a triangle with the same side-lengths lying on the plane of constant curvature equal to $c$:
$$
\alpha + \beta + \gamma \geqslant \alpha_c + \beta_c + \gamma_c.
$$
Such manifolds are called \textit{Alexandrov spaces of curvature $\geqslant c$}. For arbitrary curves in these spaces one can define a notion of an angle similar to those for the segments. We will say that a curve $\gamma$ emanating for a point $O$ has at this point a defined direction, if it makes with itself some angle, obviously, equal to zero.

Let us now introduce the notion of integral geodesic curvature. We start with integral geodesic curvature of a polygonal line composed of a shortest paths. Let $\gamma$ be such a line without self-intersections, and suppose $A_1$, $A_2$, ... , $A_n$ are the inner vertexes of $\gamma$. Given a direction on $\gamma$, we have the well-defined left and right sides. Let $\alpha_i$ be a measured from the right angle of the sector between two bars of the polygonal line meeting at the vertex $A_i$. Then the \textit{right integral geodesic curvature} of $\gamma$ is the quantity
$$
\varphi(\gamma) = \sum \limits_{i} (\pi - \alpha_i),
$$ 
where we sum over all inner vertexes. The \textit{left integral geodesic curvature} is defined in a similar way by measuring angles from the left.

Let now $\gamma$ be an arbitrary curve without self-intersections with the end points $A$ and $B$, and assume $\gamma$ has the defined directions at $A$ and $B$. After setting a direction of the curve $\gamma$, let us construct a sequence of simple polygonal lines $\gamma_n$ lying in the right semi-neighborhood of $\gamma$ and that converges to $\gamma$. Suppose $\varphi_n$ is the right geodesic curvature of the polygonal line $\gamma_n$, and $\alpha_n$, $\beta_n$ are the measure from the right angles made by the first and the last bars of the polygonal line $\gamma_n$ with the curve $\gamma$. Then the \text{right integral geodesic curvature} is defined as
$$
\varphi(\gamma) = \lim (\alpha_n + \beta_n + \varphi_n).
$$

The right geodesic curvature is defined similarly~\cite{Al48}.

A curve $\gamma$ is called \textit{$\lambda$-convex} with $\lambda>0$, if for each sub-arc $\gamma_1$ of $\gamma$
$$
\frac{\varphi(\gamma_1)}{s(\gamma_1)} \geqslant \lambda,
$$
where $s(\gamma_1)$ is the length of an arc $\gamma_1$. For regular curves in a two-dimensional Riemannian manifold this condition is equivalent to the fact that the geodesic curvature at each point $\geqslant \lambda > 0$. In the general case such a condition allows a curve to have corner points.

To prove the main theorem we will need the gluing theorem of A.\,D.~Alexandrov. Suppose that closed domains $G_1$, $G_2$, ... ,$G_n$ are cut out of manifolds of curvature $\geqslant k$ and each is bounded by a finite number of curves. We say that a manifold $R'$ with inner metric is \textit{glued from} the domains $G_1, \ldots, G_n$ if it can be partitioned into domains $G_1', \ldots, G_n'$ that are isometric to, respectively, $G_1, \ldots, G_n$.  The domains $G_1', \ldots, G_n'$ may have mutually identified parts of their boundaries; we call such parts as \textit{edges} of the domains $G_1', \ldots, G_n'$. Points were more than two domains meet together are called \textit{vertexes}. Everywhere below we assume that each domain has only finite number of edges and vertexes. The identification of edges and vertexes of $G_i'$'s can be transfered, by isometry, to the initial domains $G_i$. After that these domains by themselves generate a manifold $R$ homeomorphic to $R'$, and for which the inner metric is naturally defined. In particular, the length of a curve in the manifold $R$ is defined as the sum of the lengths of its pieces lying in each $G_i$.

A manifold glued from domains is almost completely determined by saying what parts of the boundaries should be identified. The identification rules are as follows:

\begin{enumerate}
\item
Domains, meeting at one vertex, meet in the same way as sectors, which add up to a disc, meet at the disc's center. From this assertion already follows that edges are identified pairwise.

\item
One can pass from one domain to another by going through domains with identified edges.

\item
Identified edges, as well as any identified parts of edges, have equal lengths. 
\end{enumerate}

\begin{theoremA}[A.~D.~Alexandrov, \cite{Al48}]
\label{thmA}
In order to get after gluing domains $G_1$,$\ldots$, $G_n$ with metric of curvature $\geqslant c$ a manifold $R$ with metric of curvature $\geqslant c$ it is necessary and sufficient to satisfy two conditions: 1) the sum of integral geodesic curvatures of any two identified parts of the boundaries is non-negative; 2) the sum of angles of domains meeting at a point is not bigger than $2\pi$. 
\end{theoremA}

\begin{theoremB}[A.~D.~Alexandrov, \cite{Al48}]
\label{thmB}
A metric space with inner metric of curvature $\geqslant c$ homeomorphic to a sphere is isometric to a closed convex surface in a simply connected space of constant curvature equal to $c$.
\end{theoremB}

\begin{theoremC}[A.\,V.~Pogorelov, \cite{Pog73, Pog51}]
\label{thmC}
Closed isometric convex surfaces in the three-dimensional Euclidean and spherical spaces are equal up to a rigid motion.
\end{theoremC}

\begin{theoremD}[A.\,D.~Milka, \cite{Mil80}]
\label{thmD}
Closed isometric convex surfaces in the three-dimensional Lobachevsky space are equal up to a rigid motion.
\end{theoremD}

\begin{theoremE}[W.~Meeks and S.\,T.~Yau, \cite{MYa82}]
\label{thmE}
Let $M$ be a convex three-dimensional manifold in the spherical space $S^3$, and $\gamma \subset \partial M$ be a closed Jordan curve. Then the curve $\gamma$ bounds an embedded surface which is a solution to the Plateau problem. Moreover, this surface either entirely lies on the boundary $\partial M$, or the interior of the surface lies inside $M$.
\end{theoremE}

\begin{theoremF}
Let $\gamma$ be a closed embedded $\lambda$-convex curve, with $\lambda > 0$, lying in a two-dimensional model space of constant curvature equal to $c$. If $L$ and $F$ are, respectively, the length of $\gamma$ and the area of the domain, enclosed bythe curve, then
\begin{enumerate}
\item
\emph{(A.~Borisenko, K.~Drach, \cite{BorDr14})} for the Euclidean plane, i.e. for $c = 0$,
\begin{equation*}
F \geqslant \frac{L}{2\lambda} - \frac{1}{\lambda^2} \sin \frac{L\lambda}{2};
\end{equation*}

\item
\emph{(A.~Borisenko, K.~Drach, \cite{BorDr15_1})} for the spherical space, i.e. for $c = k^2$,
\begin{equation*}
F \geqslant \frac{4}{k^2} \arctan\left(\frac{\lambda}{\sqrt{\lambda^2+k^2}}\tan\left(\frac{\sqrt{\lambda^2+k^2}}{4}L\right)\right)-\frac{L\lambda}{k^2};
\end{equation*}

\item
\emph{(K.~Drach, \cite{Dr14})} for the hyperbolic space, i.e. for $c = -k^2$, when
\begin{enumerate}
\item
$\lambda > k$,
\begin{equation*}
F \geqslant \frac{L \lambda }{k^2} - \frac{4}{k^2} \arctan\left(\frac{\lambda}{\sqrt{\lambda^2 - k^2}}\tan\left(\frac{\sqrt{\lambda^2 - k^2}}{4}L\right)\right);
\end{equation*}

\item
$\lambda = k$,
\begin{equation*}
A \geqslant   \frac{L}{k} - \frac{4}{k^2} \arctan\frac{k L}{4};
\end{equation*}

\item
$\lambda < k$,
\begin{equation*}
F \geqslant \frac{L \lambda}{k^2} - \frac{4}{k^2} \arctan\left(\frac{\lambda}{\sqrt{k^2 - \lambda^2}}\tanh \left(\frac{\sqrt{k^2 - \lambda^2}}{4}L\right)\right).
\end{equation*}
\end{enumerate}  
\end{enumerate}

\end{theoremF}

\section{Proof of the theorem}

Let us take two copies of the domain $G$. We can identify their boundary curves by isometry. After such a procedure we obtain two-dimensional manifold $F$ homeomorphic to a sphere $S^2$. By Theorem A, the constructed manifold $F$ will be a manifold with inner metric and with curvature $\geqslant c$ (in the sense of Alexandrov). By Theorem B, this manifold can be isometrically embedded as a closed convex surface $F_1$ in the simply-connected space $M^3(c)$ of constant curvature equal to $c$. From theorems C and D it follows, that up to a rigid motion this surface is unique. Moreover, the image of the curve $\gamma$ is a plane curve, and the boundary of a cap isometric to $G$. We will denote this cap and its boundary curve with the same letters.

To see that the image of $\gamma$ is a plane curve, assume the contrary. Suppose $\gamma$ is not plane. Perform the reflection of the surface $F_1$ with respect to a plane passing through 3 points on $\gamma$ that do not lie on a line. We will get the mirrored surface $F_2$ that cannot be matched with $F_1$ by a rigid motion. The domains $G_1$ and $G_2$ are mapped to domains $\tilde G_1$ and $\tilde G_2$ on $F_2$, the curve $\gamma$ is mapped to $\tilde \gamma$. But $G_1$ is isometric to $G_2$, $\tilde G_2$ is isometric to $\tilde G_1$. Let us reverse the orientation of the domains $\tilde G_1$, $\tilde G_2$. Then the surface $F_2$ will be isometric to $F_1$. By theorems C and D, $F_1$ and $F_2$ can be matched with each other with a rigid motion. And since $\gamma$ and $\tilde \gamma$ have three common points, these curves must coincide after a rigid motion. This may happen only if $\gamma$ is a plane curve, and is fixed under the reflection.
 
The cap is a graph over a plane domain $\bar G$ which is enclosed by $\gamma$. In the cases $c=0$, or $c=-k^2$ with $k>0$, by a direct computation of the area of a graph over a plane domain, we obtain that the area of $G$ is non less than the area of the plane domain $\bar G$. And the equality holds if and only if the cap $G$ coincide with $\bar G$.

By \textit{plane domains} we will understand domains on totally geodesic two-dimen\-sional surfaces in spaces of constant curvature; similarly, we will call \textit{lines} geodesic lines in these spaces.

In the spherical space $S^3(k^2)$ let us consider a domain $M$ bounded by the two-dimensional totally geodesic sphere containing $\gamma$; $M$ is a closed hemisphere. This hemisphere is a convex three-dimensional manifold in the sense of~\cite{MYa82}, thus we can apply Theorem E. By this theorem, the minimal surface $M^2$ which is a solution of the Plateau problem either coincides with the domain $\bar G$ on the boundary of $M$, or its interior lies inside the hemisphere. Since the length of the boundary curve $\gamma$ is less then $2\pi$, it lies inside a two-dimensional open totally geodesic hemisphere of the boundary $\partial M$. Assume that the minimal surface lies inside the domain $M$. This means this surface lies in some open hemisphere of the spherical space $S^3$. Consider a geodesic map this hemisphere onto euclidean space $E^3$. The curve $\gamma$ will be mapped to some plane curve $\tilde \gamma$, and the minimal surface $M^2$~-- in some regular surface $\tilde M^2$ lying in the half-space bounded by the plane of $\tilde \gamma$. From this it follows that there is a point on $\tilde M^2$ with positive Gauss curvature. Thus, by the properties of geodesic mappings, there is a point of positive extrinsic curvature on $M^2$, which is a contradiction to minimality of $M^2$.  

Therefore, the solution to the Plateau problem is the plane domain $\bar G$. From this we get that either the area of $G$ is bigger than the area of $\bar G$, or the domains $G$ and $\bar G$ coincide. Since $G$ is a cap, the integral geodesic curvature of any arc of $\gamma$ calculated on $\bar G$ is not less than the corresponding integral geodesic curvature of this calculated on the cap $G$. This means that $\gamma$ as a boundary curve of $\bar G$ is also $\lambda$-convex. 

Indeed, intrinsic curvature $\omega(D)$ of a Borel set $D$ on a convex surface in a space of constant curvature equal to $c$ is
$$
\omega(D) = \psi(D) + F(D),
$$
where $\psi(D)$ is extrinsic curvature, $F(D)$~-- area of $D$~\cite{Al48}.  Consider a closed convex surface $M^2$ bounded by $G$ and the plane domain $\bar G$, and a surface $\bar M^2$ made up from the double-covered domain $\bar G$.

Intrinsic curvature concentrated on $\gamma$ is
$$
\omega(\gamma) = \tau_\gamma (G) + \tau_\gamma (\bar G),
$$
where $\tau_\gamma (G)$ is the integral geodesic curvature of $\gamma$ calculated from the side of $G$, $\tau_\gamma (\bar G)$~-- the integral geodesic curvature calculated from the side of $\bar G$.

Since $F(\gamma) = 0$, then
\begin{equation*}
\begin{aligned}
&\psi_M(\gamma) = \tau_\gamma(G) + \tau_\gamma (\bar G),\\
&\psi_{\bar M}(\gamma) = 2\tau_\gamma (\bar G)
\end{aligned}
\end{equation*}

From the definition of the extrinsic curvature~\cite{Al48} it follows that $\psi_{\bar M} (\gamma) \geqslant \psi_M(\gamma)$, since each supporting to $M^2$ plane at a point on $\gamma$ is also supporting to $\bar M^2$. Thus we obtain $\tau_\gamma(\bar G) \geqslant \tau_\gamma (G)$. Moreover, this inequality holds for any sub-arc of $\gamma$ as well. For the Euclidean space this inequality was proved by V.\,A.~Zalgaller~\cite{Za52}. 

Consider the domain $\bar G$ with the boundary curve $\gamma$. This domain lies in a two-dimensional space of constant curvature $c$ and satisfies the conditions of the theorem. Applying Theorem F together with the considerations above, we get the proof of the theorem.


\end{document}